\documentclass[11pt,a4paper]{article}
\usepackage{epsfig,amsmath,amssymb,enumerate,mathrsfs,theorem}
\newtheorem{thm}{Theorem}
\newtheorem{cor}[thm]{Corollary}
\newtheorem{lem}[thm]{Lemma}

\def\pn{\par\noindent}

\title{Fermat's Last Theorem for the Exponent $3$}
\author{Roy Barbara }
\date{}
\begin{document}
\maketitle
$ $\\
In this note we present a new proof (in Euler's style) of F.L.T. in the case $n=3$: The crucial lemma receives a short proof and the classical  "descent" in Two main cases boils down to a One-case proof.\\ \\
\textbf{\hspace*{.3cm}Notations}: $\mathbb{Z}^*$ stands for $\mathbb{Z}-\{0\}$ and $\mathbb{N}^*$ for $\mathbb{N}-\{0\}$.\, Throughout, $a,\,b,\,p,\,q\,$ lie in $\mathbb{Z}^*$. \, \hspace*{.8cm} $(a,\,b)$ denotes the g.\,c.\,d.
\begin{enumerate}
\item \textbf{Introduction:}\\
  Recall that every prime $p=6k+1$ of $\mathbb{N}$ has a unique representation as $p=r^2+3s^2,\;r,\,s\,\in\,\mathbb{N}^*$. Set $\displaystyle \omega = \frac{-1+\sqrt{-3}}{2},\,A=\mathbb{Z}(\omega)=\{\displaystyle \frac{m+n\sqrt{-3}}{2},\,m,\,n\,\in\,\mathbb{Z},\,m+n\,\,even\}\,\supset\,\mathbb{Z}(\sqrt{-3})$. The units of $A$ are $\pm1,\,\pm\omega,\,\pm\omega^2$. A is a euclidean (hence factorial) ring for the norm $N(\alpha+\beta \sqrt{-3})=\alpha^2+3\beta^2$. Define the "\textit{atoms}" of $A$, \textit{all lying} in $\mathbb{Z}(\sqrt{-3}),$ as the elements: \,$2$,\, \,$p=6k-1$ prime of $\mathbb{N},\,\sqrt{-3},\,r+s\sqrt{-3}$ and $r-s\sqrt{-3}$ for each prime $p=6k+1$ of $\mathbb{N},\,p=r^2+3s^2,\,r,\,s\in\,\mathbb{N}^*$.\,
  The atoms of $A$,\, pairwise coprime, represent (up to units) \textit{all} the irreductible (or prime) elements of $A$. \,
  Let $z\,\in\,A-\{0\},\,z$ not a unit. Then, $z$ has a unique "\textit{atomic}" decomposition as $z=\epsilon . \pi$, where $\epsilon$ is a unit of $A$ and $\pi\in \mathbb{Z}(\sqrt{-3})$ is a product of atoms in a \textit{unique} way (up to order). Every  atom of $z$ corresponds to \textit{one} prime factor (or its square) of $N(z):\; $ $2$ \,\textit{to} \,$2^2$, \,$p=6k-1$  \textit{to}  $p^2$, \, $\sqrt{-3}$ \, \textit{to} \,$3$, \, $r+s\sqrt{-3}$, \textit{resp.} $r-s\sqrt{-3}$, \,\textit{to} \, $p\, (=6k+1)=r^2+3s^2$.\\
  The following properties are trivial:
  \begin{enumerate}[(R1)]
     \item If $a+b$ is odd, then, the atom $2$ does not divide $a+b\sqrt{-3}$ in $A$.
     \item If $(a,b)=1$, then, the atom $p=6k-1$ does not divide $a+b\sqrt{-3}$ in $A$.
     \item if $(a,b)=1$, then, the square of the atom $\sqrt{-3}$\, does not divide $a+b\sqrt{-3}$ in $A$.
  \end{enumerate}
  In particular we have:
  \begin{enumerate}[(P1)]
  \item If $a+b$ is old and $(a,b)=1$, then, $\textit{N}(a+b\displaystyle\sqrt{-3})=a^2+3b^2$ \, has neither factor $2$, \, nor $p=6k-1$, \, nor $3^2$.
      \end{enumerate}
\item \textbf{Proof of F.L.T. in the case $n=3$}:\\
  \textit{Notice} the following properties:
  \begin{enumerate}[(P2)]
    \item If $(a,b)=1$ and $p=6k+1=r^2+3s^2$\, is a prime of $\mathbb{N}$, then, the product of the atoms $\lambda = r+s\sqrt{-3}$ and $\overline{\lambda}=r-s\sqrt{-3}$ \, does not divide $z=a+b\sqrt{-3}$ \,in $A$.\\
        \textbf{\textit{Proof}}: Otherwise, for some $m,\,n\,\in \mathbb{Z}: \; a+b\sqrt{-3}=\lambda \overline{\lambda}\displaystyle (\frac{m+n\sqrt{-3}}{2})=p\displaystyle (\frac{m+n\sqrt{-3}}{2})$.\\
        Then, $2a=pm,\; 2b=pn$. \textit{Hence}, $p/a$ and $p/b$, a contradiction.
     \end{enumerate}
  \begin{enumerate}[(P3)]
    \item Let $a+b$ odd and $(a,\,b)=1$. \,Suppose that $a^2+3b^2=c^3,\;\,c\in\mathbb{N}^*$.  Then, $z=a+b\sqrt{-3}$\; is a cube in $\mathbb{Z}(\sqrt{-3})$.\\
       \textbf{\textit{Proof}}: Let $z=\epsilon .\pi$ be the atomic decomposition. As $N(z)=c^3$, by $P1,\,\, c$ has only prime factors $p=6k+1$: For \textit{each} such $p,\,\, p.p.p $ divides $N(z)$. By P2, \textit{no two} of the $3$ corresponding atoms of $z$ can be \textit{conjugate}. Hence, these 3 atoms are equal. This argument shows that $\pi$ is a cube in $\mathbb{Z}(\sqrt{-3})$. Hence, $z=\epsilon .\gamma ^3,\, \gamma\, \in\, \mathbb{Z}(\sqrt{-3})$. Then, $\epsilon^{-1}z=\gamma^3\,\in\,\mathbb{Z}(\sqrt{-3})$. But, since $a+b$ is odd, $\pm \,\omega z, \,\pm\, \omega^2z\,\notin \,\mathbb{Z}(\sqrt{-3})$. \, Hence $\epsilon^{-1}=\pm1$, so $z=(\pm\gamma)^3$.
  \end{enumerate}
  $ $\\
       \begin{lem} Let $a+b$ odd and $(a,\,b)=1$. Suppose that $a^2+3b^2=c^3,\,c\,\in \mathbb{N}^*$. \,Then,
       \begin{enumerate}[(i)]
         \item $3$ divides $b$.
         \item For some $e,\,f\,\in \mathbb{Z}^*\;:\, a=e(e^2-9f^2)$\, and \, $b=3f(e^2-f^2)$.
       \end{enumerate}
       \end{lem}
       \textit{\textbf{Proof:}}\; By P3, for some $e+f\sqrt{-3}\,\in\mathbb{Z}(\sqrt{-3}),\, a+b\sqrt{-3}=(e+f\sqrt{-3})^3$.\\
       Clearly, $e,\, f\neq\,0.\, \,\,$ By expanding we get (ii). \hspace{.2cm} (i) follows.
  \begin{enumerate}[(P4)]
    \item Let $p+q$ odd and $(p,\,q)=1$. Suppose that $2p(p^2+3q^2)=\theta^3,\,\theta\,\in\mathbb{Z}^*$. \,Then, $q$ is odd and $9 \mid p.\,$ \textit{In particular}, $3\,\mid  \,\theta $ and $3\nmid q$.\\ \\
        \textbf{\textit{Proof}:}\\$8\mid \theta^3$\, so $4\mid p(p^2+3q^2)$. As $p^2+3q^2$ is odd, then $4\mid p$. Hence $q$ is odd.

        If $3\nmid p$, then $3\nmid \theta, $ so $\theta^3 \equiv\,\pm\,1\,(mod.9)$ and $2p^3\,\equiv \,\pm\,2\,(mod.9)$. Now, $(2p,\,p^2+3q^2)=1$ so $p^2+3q^2$ is a cube in $\mathbb{N}^*$. By the lemma $3\mid q,$ so $2p(p^2+3q^2)\equiv 2p^3\equiv \pm2(mod.9)$, a contradiction.\\
        Hence $3\mid p$, so $3\mid \theta$ and $3\nmid q$. Now $27\mid p(p^2+3q^2)$ and $9\nmid p^2+3q^2$ \,\,yield $9\mid p$.\\
  \end{enumerate}
        \begin{thm}
          The equation $2p(p^2+3q^2)=\theta^3$\,\,\, (1)\,\,\, with $p+q$ odd and $(p,\,q)=1$ is insolvable in $\mathbb{Z}^*$.
        \end{thm}
        \textbf{Proof (by descent):} Consider a solution as in \,(1)\, with $\mid \theta \mid$ minimum. By P4,\, $q$ is odd,\, $9\mid p,\,\,3\mid \theta,\,\, 3\nmid q$. We may write $$\displaystyle \big{(}\frac{2p}{9}\big{)}\big{(}\frac{p^2+3q^2}{3}\big{)}=\big{(}\frac{\theta}{3}\big{)}^3\;\;\;\;\;\;\;\;\;\;\;\;(2)$$
        As $3\nmid \displaystyle \frac{p^2+3q^2}{3}$, clearly, $\displaystyle \big{(}\frac{2p}{9},\,\frac{p^2+3q^2}{3}\big{)}=1.$ \, Hence, $$\displaystyle  \frac{2p}{9}=u^3\;\;\;\;\;(3)\;\;\;\;\; \frac{p^2+3q^2}{3}=v^3\;\;\;\;\; (4)\;\;\;\;\; uv=\frac{\theta}{3}\,\,\, (u,\,v\,\in \mathbb{Z}^*)$$
        From \, $q^2+3\displaystyle (\frac{p}{3})^2=v^3\,\,\,(4)\,\,\,$ and the lemma, we get $$q=e(e^2-9f^2)\;\;\;\;\;\;\;\; (5)
         \;\;\;\;\;\;\;\;\displaystyle \frac{p}{3}=3f(e^2-f^2)\;\;\;\;\;\;\;\; (6)$$\\
         Since $(q,\,\frac{p}{3})=1$, by (5) and (6), $(e,f)=1.\,$ As $q$ is odd, by (5), $e$ is odd and $f$ even. By (6) and (3) we have $\displaystyle \frac{2p}{9}=(e+f)(e-f)(2f)=u^3$, where the $3$ factors are pairwise coprime. Hence, $e+f=r^3,\;\,e-f=s^3,\,\; 2f=t^3,\;\; rst=u,\;\;\; r,\,s$ odd,\, $(r,\,s)=1$.
         Set $\alpha=\displaystyle \frac{r-s}{2},\, \beta=\frac{r+s}{2}.\;\,\,\,\alpha,\,\beta\neq 0,\, \alpha+\beta$ is odd and $(\alpha,\,\beta)=1.$ \,Rewrite $r^3-s^3=t^3$ as $$2\alpha(\alpha^2+3\beta^2)=t^3$$
         Where $\mid t\mid \,<\, \mid \theta \mid$  since $rstv=uv=\displaystyle \frac{\theta}{3}$.\\
         \begin{cor}
             The cubic Fermat's equation $x^3+y^3+z^3=0$ is insolvable in $\mathbb{Z}^*.$
         \end{cor}
            \textbf{\textit{Proof}:} Otherwise, pick a primitive solution with $x,\,y$ odd. Clearly, $x\neq\,y.$  Putting $\displaystyle p=\frac{x+y}{2}$ and $\displaystyle q=\frac{x-y}{2}$ \;yields\; $2p(p^2+3q^2)=(-z)^3$,\; where $p+q$ is odd and $(p,\,q)=1$.
\item \textbf{Some Applications.}\\
Applications to F.L.T. in the case $n=3$ are numerous: We give some illustrative examples.
\begin{enumerate}[(i)]
\item Let $p,\,q\,\in \mathbb{Z}$ with $p\neq\,0$ and $q\neq p^2.\,$ Then, $12q^3-3p^6$ is not a perfect square in $\mathbb{Z}$.\\
\textit{Indeed}, for the purpose of contradiction, suppose that $12q^3-3p^6=r^2,\,r\in \mathbb{Z}$. Clearly, $q\neq0$. Set $x=3p^3+r,\,\,y=3p^3-r,\,\,z=6pq,$ and note that $x,\,y,\,z$ lie in $\mathbb{Z}^*$. Then, using the identity $(u+v)^3+(u-v)^3=2u^3+6uv^2$, we obtain $x^3+y^3=z^3$, a contadiction.
\item The diophantine equation $9x^3=zy^3+z^2+6xyz$ has no solution in $\mathbb{Z}^*$.\\
\textit{Indeed}, suppose this equation holds for some non-zero integers $x,\,y,\,z$. Set $p=y,\,q=y^2+3x,\,r=3(y^3+6xy+2z)$. Then, $p\neq 0$ and $q\neq p^2$. After checking carefully the relation $12q^3-3p^6=r^2$, we obtain a contradiction in virtue of the property in (i).
\item The equation $\big{(} \frac{1}{2}(a+b+c)\big{)}^3=3abc$ has no solution in non-zero integers.\\
\textit{Indeed}, suppose that $a,\,b,\,c\,\in \mathbb{Z}^*$ satisfy\,\,\,\,$(a+b+c)^3=24abc$.\hspace*{0.5 cm} \hspace*{2.2cm}(7)\\ Expanding (7), we get \begin{center} $a^3+b^3+c^3-18abc+3a^2(b+c)+3b^2(c+a)+3c^2(a+b)=0$ \hspace{2.1cm}\,\;\;\;(8)\end{center}
Set $r=a+b-c,\,s=a-b+c,\,\,t=-a+b+c$. One checks that $r\neq 0$ ($r=0$ leads by $(7)$ to $a^2+b^2-ab=0)$ and similarly that $s\neq 0$ and $t\neq 0$. Now, expanding $r^3+s^3+t^3=(a+b-c)^3+(a-b+c)^3+(-a+b+c)^3$, we obtain \begin{center}$r^3+s^3+t^3=a^3+b^3+c^3-18abc+3a^2(b+c)+3b^2(c+a)+3c^2(a+b)\;\;\;\;\; \;\;\;(9)$\end{center}
Comparing (8) and (9) provides $r^3+s^3+t^3=0$, a contradiction.
\end{enumerate}
\end{enumerate}

\bigskip
{\footnotesize \pn{\bf Roy Barbara. } \;
\\Lebanese University,
Faculty of Science II.\\
Fanar Campus. P.O.Box 90656. \\
Jdeidet El Metn. Lebanon.\\
{\tt Email: roybarbara.math@gmail.com \\


\begin{thebibliography}{99}
%
\bibitem{2} Ribenboim, P.,  \textit{13 lectures on fermat's Last Theorem}. Springer-Verlag, 1999.
\bibitem{2} Legendre, A.M., \textit{Th\'eorie des nombres}, vol. 2, reprint, blanchard, Paris 1955.
\bibitem{2} Edwards. H.M. \textit{Fermat's Last Theorem: A Genetic Introduction to Algebraic\\ Number Theory.} New York: Springer-Verlag, 1977.
\bibitem{2} Mordell, L.J. \textit{Three Lectures on Fermat's Last Theorem.} New York: Chelsea, 1956.
\bibitem{2} Cox, D.A.\hspace{.3cm} "Introduction to Fermat's Last Theorem".\hspace{.3cm} \textit{Amer. Math. Monthly} 101, \\ 3-14, 1994.
\bibitem{2} Ribenboim, P. \textit{Fermat's Last Theorem for Amateurs}. New York: Springer-Verlag, 1999.
\bibitem{2} I. Niven, H. Zuckerman. \textit{An Introduction to the Theory of Numbers}, $4^{th}$ Ed., John Wiley \& sons, Inc., New York, 1980.
\end{thebibliography}
\end{document}